%% file: CDC_Jasour.tex
\documentclass[letterpaper, 10 pt, conference]{ieeeconf}  

\IEEEoverridecommandlockouts                              
\usepackage{fullpage,enumerate,setspace,changepage,multirow,rotating,xcolor}
\usepackage{graphicx,amsmath,amssymb,amsfonts,url, algorithm, algorithmic, etoolbox}
\usepackage{array,subfigure}
\usepackage[small]{caption}
\AtBeginEnvironment{align}{\setcounter{subeqn}{0}}
\newcounter{subeqn} \renewcommand{\thesubeqn}{\theequation\alph{subeqn}}%
\newcommand{\subeqn}{%
	\refstepcounter{subeqn}
	\tag{\thesubeqn}
}
\input defs.tex

\newtheorem{lemma}{Lemma}

\def\norm#1{\|#1\|}

\def\norm#1{\|#1\|}

\makeatletter
\newcommand{\thickhline}{%
	\noalign {\ifnum 0=`}\fi \hrule height 1pt
	\futurelet \reserved@a \@xhline
}

\title{\LARGE \bf
Moment-Sum-Of-Squares Approach For Fast Risk Estimation In Uncertain Environments
}

\author{Ashkan M. Jasour, Andreas Hofmann, Brian C. Williams\\
	MIT, Computer Science and Artificial Intelligence Laboratory \\
	\{jasour,hofma,williams@csail.mit.edu\}
\thanks{
	This work was supported in part by Boeing grant MIT-BA-GTA-1 and by the Toyota Research Institute (TRI). However, this article solely reflects the opinions and conclusions of its authors and not TRI or any other Toyota entity.}
}

\begin{document}

\onecolumn

\maketitle
\thispagestyle{empty}
\pagestyle{empty}


\begin{abstract}
	In this paper, we address the risk estimation
	problem where one aims at estimating the probability of
	violation of safety constraints for a robot in the presence of bounded
	uncertainties with arbitrary probability distributions. In this problem, an unsafe set is described by level sets of polynomials that is, in general, a non-convex set. Uncertainty arises due to the probabilistic parameters of the unsafe set and probabilistic states of the robot. To solve this problem, we use a moment-based representation of probability distributions. We describe upper and lower bounds of the risk in terms of a linear weighted sum of the moments. Weights are coefficients of a univariate Chebyshev polynomial obtained by solving a sum-of-squares optimization problem in the offline step. Hence, given a finite number of moments of probability distributions, risk can be estimated in real-time. We demonstrate the performance of the provided approach by solving probabilistic collision checking problems where we aim to find the probability of collision of a robot with a non-convex obstacle in the presence of probabilistic uncertainties in the location of the robot and size, location, and geometry of the obstacle.
\end{abstract}
\section{Introduction}
Probabilistic analysis plays a key role in planning and control problems in the presence of uncertainties.
In this paper, we consider the risk estimation problem where we aim to estimate the probability of failure in the presence of uncertainties. 
This problem has many applications in different areas. For example, in probabilistic motion planning of robots where we need the probability of collision with obstacles in uncertain environments (\cite{Bool1,Bool2,Bool3,Bool4,Bool5}). Another application is in stochastic scheduling problems where the probability of feasibility of each plan given the probabilistic temporal constraints is required (\cite{sch1,sch2,sch3}).\\ 

The problem of computing the probability is challenging because it requires evaluation of multivariate integrals over non-convex sets. Several approaches have been proposed to find the probability over a given set. However, the proposed approaches are limited to particular uncertainties and sets. For example, Boole's inequality is widely used to estimate the probability of violation of linear constraints (\cite{Bool1,Bool2,Bool3}). More precisely, the probability of being safe in the presence of a convex polytopic obstacle $\chi$ represented by a conjunction of linear inequality constraints, e.g., $\chi=\{\cap_{j=1}^N \chi_j\}$, $\chi_j=\{ x\in \mathbb{R}^n : a'_jx \leq b_j \}$, is calculated as follows: $ 1-\mbox{probability} ( \cap_{j=1}^{N} \chi_j ) = \mbox{probability} ( \cup_{j=1}^{N} \bar{\chi}_j ) \leq \sum_{j=1}^N \mbox{probability} (\bar{\chi}_j)$ where $\bar{\chi}_j=\mathbb{R}^n \setminus \chi_j$ is the complement set. This results in a conservative upper bound on the probability. Uncertainty sampling based methods are also widely used in many probabilistic planning and control applications (\cite{Samp1,Samp2,Samp3}). Being a randomized approach, no analytical bounds can be provided on the probability.  
In (\cite{ Opt1, Ind1, Ashk}), semidefinite programs (SDP) are provided to estimate the probability of semialgebraic sets. These methods rely on polynomial approximation of an $n$-dimensional indicator function of a given set that are formulated as a sum of squares (SOS) optimization. The obtained SDPs easily become computationally intractable as the size of the original problem increases.\\


In this paper, we leverage SOS based techniques to provide upper and lower bounds of the probability of violation of safety constraints described by level sets of $n$-variate polynomials. The proposed method can deal with bounded uncertainties with arbitrary probability distributions and also uncertain nonconvex safety constraints e.g., obstacles with uncertain location, size, and geometry. The provided method relies on a convex optimization that looks for a univariate polynomial indicator function. Using the proposed approach, we describe upper and lower bounds of the risk as a linear weighted sum of the moments of uncertainties. The weights are coefficients of a univariate Chebyshev polynomial obtained by solving a univariate SOS optimization. 

The key innovations of our proposed approach are as follow: (1)  the proposed approach performs
numerical computations in the offline step, 
and uses these results to efficiently compute the risk bounds for
the given moments of probability distributions, in real-time. Hence, in the presence of time varying or state-dependent uncertainties, it can update the risk bounds by only updating the moment information, in real-time, 
(2) to reduce the size of SOS optimization, the proposed approach solves a univariate SOS optimization. Hence, computation time reduces significantly compared to the multivariate SOS based techniques.\\

The outline of the paper is as follows: in Section 2, we cover the notation adopted in the paper, and present preliminary results on polynomials; Section 3 includes the problem statement and a motivating example; Section 4 details the proposed technique to estimate the probability with an illustrative example; in Section 5, we present numerical results, followed by some concluding remarks given in Section 6.


\section{ Notation and Preliminary Results}\label{Notation}

\label{sec:definitions}

This section covers notation and includes some basic definitions of polynomials and moments (\cite{SOS1,SOS2,SOS3,Ind1,Ashk}). Given $n$ and $d$ in $\mathbb{N}$, we define {\small $S_{n,d} := \binom{d+n}{n}$} and {\small $\mathbb{N} ^{\rm n}_d := \{\alpha \in \mathbb N^n : \norm{\alpha}_1 \leq d \}$}. Also, given two sets $A$ and $B$, we define the set difference by {\small $A \setminus B = \{x: x\in A, x \notin B\}$}.\\

\textbf{Standard Polynomials:} Let {\small $\mathbb{R}[x]$} be the set of real polynomials in the variables {\small $x \in \mathbb{R}^n$}. Given {\small $\cP\in\mathbb{R}[x]$}, we represent {\small $\cP$} as {\small $\sum_{\alpha\in\mathbb{N}^n} p_\alpha x^\alpha$} using the standard basis {\small $\{x^\alpha\}_{\alpha\in \mathbb{N}^n}$} of $\mathbb{R}[x]$, and {\small $\mathbf{p}=\{p_\alpha\}_{\alpha\in\mathbb{N}^n}$} denotes the polynomial coefficients. Also, let $\mathbb R_{\rm d}[x] \subset \mathbb R [x]$ denotes the set of polynomials of degree at most $d\in \mathbb{N}$. Any given {\small $\cP\in\mathbb R_{\rm d}[x]$}, has {\small $S_{n,d}$} number of coefficients.
\\

\textbf{Chebyshev Polynomials:} Chebyshev polynomials of the first kind with degree $d$ are defined as {\small $T_d(x)=cos(d \cos^{-1}(x)), \ x\in[-1,\ 1],d \in \mathbb{N}$}, \cite{cheb}. Chebyshev polynomial $T_d(x)$ can be represented in terms of powers of $x$ as {\footnotesize $T_d(x)= \frac{d}{2}\sum_{i=0}^{[d/2]} (-1)^i \frac{(d-i-1)!}{i!(d-2i)!}(2x)^{d-2i}$} (e.g., {\footnotesize$T_0(x)=1,T_1(x)=x, T_2(x)=2x^2-1$}). Also, the product of Chebyshev polynomials can be expanded as follows: {\footnotesize$T_{d_1}(x)T_{d_2}(x)=\frac{1}{2}(T_{d_1+d_2}+T_{|d_1-d_2|})$}. The important property of Chebyshev polynomials is orthogonality. 
\\

\textbf{Sum of Squares Polynomials:} Let {\small$\mathbb{S}^2[x] \subset \mathbb R [x]$} be the set of sum of squares (SOS) polynomials. Polynomial  {\small$s:\reals^n\rightarrow\reals$} is an SOS polynomial if it can be written as a sum of \emph{finitely} many squared polynomials, i.e., {\small$s(x)= \sum_{j=1}^{\ell} h_j(x)^2$} for some $\ell<\infty$ and $h_j\in\reals[x]$ for $1\leq j\leq \ell$. The following lemma gives a sufficient condition for $\mathcal{P} \in \mathbb{R}[x]$ to be nonnegative on the compact set {\small $\mathcal{K}=\{ x \in \reals^n: \mathcal{P}_j(x) \geq 0, j=1,2,...,\ell \}$}, where {\small $\mathcal{P}_j \in \mathbb{R}[x]$} (\cite{SOS1,SOS2,SOS3,Ind1}).
\begin{lemma}
	\label{sec2:lem7}
	If {\small $\mathcal{P} \in \mathbb{R}[x]$} is strictly positive on $\cK$, then {\small $\mathcal{P}$} has the SOS representation as follows:
	\begin{equation*}
	\mathcal{P} = s_0 + \sum_{j=1}^{\ell} s_j \mathcal{P}_j, \ s_j \in \mathbb{S}^2[x], \ j=0,...,\ell
	\end{equation*}
\end{lemma}
The SOS condition is a convex constraint that can be represented as a linear matrix inequality in terms of coefficients of polynomial {\small $\mathcal{P}$}.\\


\textbf{Moments of Probability Distributions} Let {\small $x=[x_1,...,x_n] \in \mathbb{R} ^n$} be a multivariate random variable with probability distribution {\small$\mu_x$}. Support of the probability distribution $\mu_x$ is denoted by $supp(\mu_x)$, i.e., the smallest closed set that contains all the sets with nonzero probability. Given {\small $\overrightarrow{\alpha}=(\alpha_1,...,\alpha_n)$} with {\small$\alpha_i$} in $\mathbb{N}$, the moment of order {\small$\overrightarrow{\alpha}$} of {\small$\mu_x$} is defined as {\footnotesize$m^{x}_{\alpha_1,\alpha_2,...,\alpha_n}=E[x_1^{\alpha_1}x_2^{\alpha_2}...x_n^{\alpha_n}]=\int x_1^{\alpha_1}x_2^{\alpha_2}...x_n^{\alpha_n} \mu_x dx_1...dx_n$}. If {\small$\mu_x$} is defined on the hyper-cube {\small$[-1, 1]^n$}, then its moments are bounded in {\small$[-1,1]$} (\cite{SOS1,Ashk}). Moments of {\small$\mu_x$} can be written in terms of the Chebyshev basis as {\footnotesize$m^{x}_{T_{\alpha_1,...,\alpha_n}}=E[T_{\alpha_1}...T_{\alpha_n}]$}. Using the mapping between Chebyshev and standard polynomials, moments in the Chebyshev basis  can be written in terms of the moments in the standard basis (e.g., {\footnotesize$n=1$, $m^x_{T_2}=2m^x_2-m^x_0$} ) \cite{Ashk}.


\section{Problem Statement}\label{Formu}

In this paper, we consider the \emph{risk estimation problem} defined as follows: let $x\in \reals^n$ be a multivariate random variable with known probability distribution $\mu_x$ defined on a compact set (e.g., uncertain position of a robot in work/joint space). The uncertain unsafe set $\chi$ (e.g., obstacle with uncertain location/size/geometry) is defined as level-sets of polynomials as follows:
\begin{equation} \label{intro_set1}
\chi(q):= \left\lbrace  x\in \reals^n: l_j \leq \mathcal{P}_{j}(x,q)\leq u_j , j=1,\dots ,\ell \ \right\rbrace 
\end{equation}
where $\mathcal{P}_{j}:{{\mathbb{R}}}^n\times {{\mathbb{R}}}^m\rightarrow{\mathbb{R}}$, $j=1,2,\dots ,\ell$ are given polynomials, $q \in \reals^m$ is a multivariate random variable with known probability distribution $\mu_q$ defined on a compact set, and $l_j, u_j \in \reals$ for $j=1,...,\ell$.
Set $\chi(q)$ is in general a non-convex set. Given the probability distributions $\mu_x$ and $\mu_q$ and unsafe set $\chi$, we focus on solving the following problem:
\begin{align} \label{intro_P1}
\mathbf{P_{risk}^*} := \mbox{Probability}_{\mu_x,\mu_q} \{x \in \chi(q)\}
\end{align}
where $\mathbf{P_{risk}^*}$ is the probability of failure due to violation of the safety constraints (e.g., probability of collision with the obstacle). The probability in \eqref{intro_P1} involves a multivariate integral over a nonconvex set i.e., {\small$\int _{\left\lbrace  (x,q): \{l_j \leq \mathcal{P}_{j}(x,q)\leq u_j\}_{j=1}^{\ell}  \ \right\rbrace } \mu_x\mu_qdx_1...dx_ndq_1...dq_m$}, which is computationally challenging.\\

\begin{figure}
	\centering
	\includegraphics[width=10cm, height=4cm]{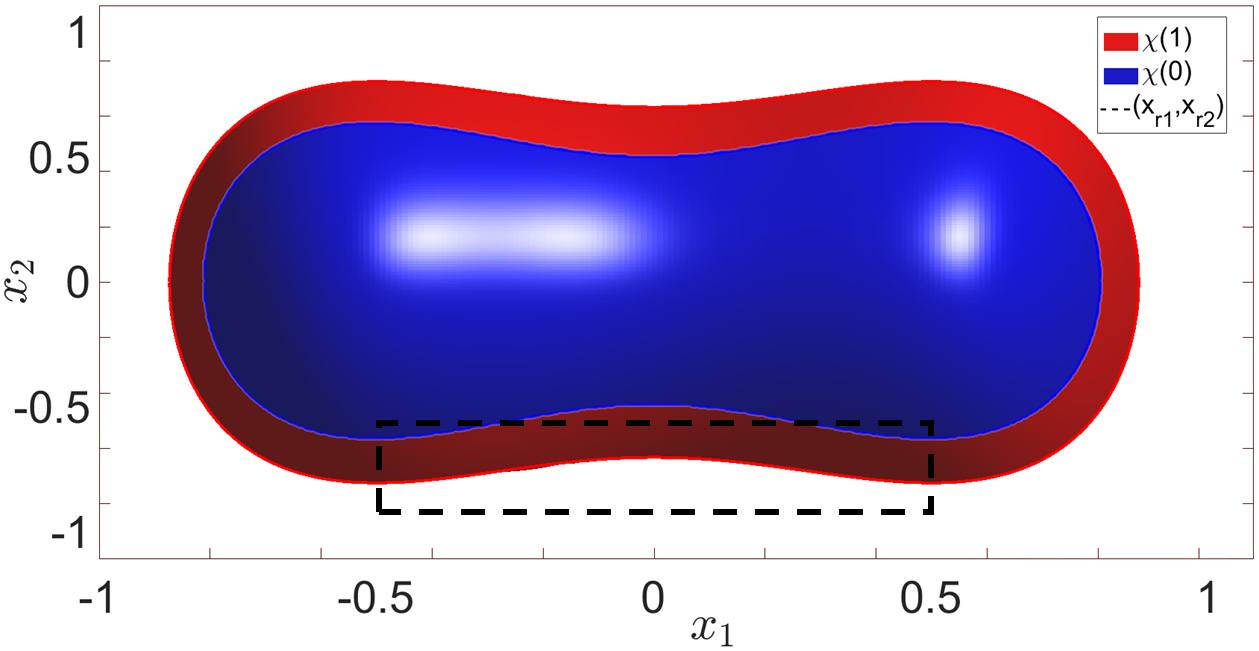}
	\caption{\small{Uncertain obstacle for parameter $q=0$ (blue) and $q=1$ (red) and possible locations of the rover (dashed line)}}
	\label{fig:E1}
\end{figure}

\textbf{Motivating Example:}
An uncertain non-convex obstacle, shown in Figure \ref{fig:E1}, is described as
{\small $\chi(q)=\{x\in\reals^2: -0.1\leq-x_1^4+0.5(x_1^2-x_2^2)+0.1q \leq 0.2\}$} where $q$ has a $Beta(4,4)$ probability distribution defined on $[0,1]$. A rover is located at $(x_{r1},x_{r2})$ where $x_{r1}$ and $x_{r2}$ have uniform probability distributions on $[-0.5,0.5]$ and $[-0.8,-0.5]$, respectively. We want to find the risk defined as the probability of collision with the obstacle, i.e., {\small $\mathbf{P_{risk}^*} = \mbox{Probability}_{\mu_{x_{r1}},\mu_{x_{r2}},\mu_q} \{-0.1\leq-x_{r1}^4+0.5(x_{r1}^2-x_{r2}^2)+0.1q \leq 0.2\}$}.



\section{Moment-Sum-Of-Squares Formulation}\label{sec:Primal}

To solve the risk estimation problem defined in \eqref{intro_P1}, we first provide a sum of squares (SOS) optimization approach involving multivariate polynomial approximation. Then, we reduce the size of the optimization problem and look for a univariate polynomial to obtain the solution of the original problem in \eqref{intro_P1}.\\

Given the polynomials $\mathcal{P}_{j}(x,q)$, $j=1,2,\dots ,\ell$, of the set in \eqref{intro_set1}, we define:
	\begin{equation} \label{eq: setK}
	\mathcal{K}= \left\lbrace  (x,q) \in \reals^n \times \reals^m, l_j \leq \mathcal{P}_{j}(x,q)\leq u_j , j=1,\dots ,\ell \ \right\rbrace 
	\end{equation}
\textit{Assumption} 1. Set $\mathcal{K}$ is a compact set. 
Hence, the projection of $\mathcal{K}$ onto $x$-coordinates denoted by $\Pi_x$ and onto $q$-coordinates denoted by $\Pi_q$ are also compact. Therefore, after rescaling of the polynomials, we assume without loss of generality that $\Pi_x \subset \textbf{X}=[-1,1]^n$ and  $\Pi_q \subset \textbf{Q}=[-1,1]^m$; Hence, $\mathcal{K} \subset \mathcal{B}=\textbf{X}\times\textbf{Q}=[-1,1]^n \times [-1,1]^m $, \cite{Ashk}. 

\textit{Assumption} 2. We assume that $supp(\mu_x) \subset \textbf{X}$ and $supp(\mu_q) \subset \textbf{Q}$. Also, moments of any order of the probability distributions can be computed \cite{Ashk}.

\textit{Assumption} 3. After rescaling of the polynomials, we assume without loss of generality that the polynomials of set $\mathcal{K}$ are bounded as $-1 \leq \mathcal{P}_{j}(x,q) \leq 1,j=1,...,\ell$ on $\mathcal{B}$; Hence, $-1 \leq l_j,u_j \leq -1, j=1,...,\ell$.

\subsection{Moment-SOS Based Risk Bounds}

Consider the defined set $\mathcal{K}$ in \eqref{eq: setK}. Let $\mathcal{I}_\mathcal{K}$ be the indicator function of the set $\mathcal{K}$, (e.g., $\mathcal{I}_\mathcal{K}=1 \ \forall (x,q)\in \mathcal{K}, \ \mathcal{I}_\mathcal{K}=0 \ \forall (x,q) \notin \mathcal{K} $). Then, the probability in \eqref{intro_P1} can be written as the following expectation:
\begin{equation}\label{eq: expect}
\mathbf{P_{risk}^*} = E[\mathcal{I}_\mathcal{K}]= \int \mathcal{I}_\mathcal{K} \mu_x\mu_qdx_1...dx_ndq_1...dq_m
\end{equation}
To evaluate the integral in \eqref{eq: expect}, one can use the polynomial approximation of the indicator function denoted by $\mathcal{P}_{\mathcal{K}}(x,q)$, as follows (\cite{Ind1},\cite{Ind2} Lemma 1): 
	\begin{align}\label{Ind}
	& \mbox{min}_{\mathcal{P}_{\mathcal{K}}(x,q)\in \mathbb{R}_d[x,q]}  \int_{\mathcal{B}} \mathcal{P}_{\mathcal{K}}(x,q)dx_1...dx_ndq_1...dq_m \\
	\mbox{s.t} &  \ \ \mathcal{P}_{\mathcal{K}}(x,q) \geq 1 \ \ \mbox{on} \ \ \mathcal{K} \subeqn \label{Ind_cons1} \\
	& \ \ \mathcal{P}_{\mathcal{K}}(x,q) \geq 0 \ \ \mbox{on} \ \ \mathcal{B} \subeqn \label{Ind_cons2}
	\end{align}	
where $\mathcal{B}$ is the bounding box defined in \textit{Assumption} 1. Note that the constraints in \eqref{Ind_cons1} and \eqref{Ind_cons2} are polynomial non-negativity constraints that can be formulated as SOS convex constraints. Similarly, one can find the polynomial approximation of the indicator function of the complement set $\bar{\mathcal{K}}= \mathcal{B} \setminus \mathcal{K} $ denoted by $\mathcal{P}_{\bar{\mathcal{K}}}(x,q)$. Obtaining degree-$d$ polynomials $\mathcal{P}_{\mathcal{K}}(x,q)$ and $\mathcal{P}_{\bar{\mathcal{K}}}(x,q) $, the following results hold:\\

\begin{lemma}\label{Lemma2}
	The risk defined in \eqref{intro_P1} is bounded by
	\begin{equation} \label{risk_bound}
	1- E[\mathcal{P}_{\bar{\mathcal{K}}}(x,q)] \leq \mathbf{P_{risk}^*} \leq E[\mathcal{P}_{\mathcal{K}}(x,q)]
	\end{equation}
	and {\footnotesize $\lim_{d \rightarrow \infty} E[\mathcal{P}_{\mathcal{K}}(x,q)] = \lim_{d \rightarrow \infty} 1- E[\mathcal{P}_{\bar{\mathcal{K}}}(x,q)] = \mathbf{P_{risk}^*}$}.
\end{lemma}

\textit{Sketch of the proof:}
$\mathcal{P}_{\mathcal{K}}(x,q) \in \mathbb{R}_d[x,q] $ is an upper bound approximation of $\mathcal{I}_{\mathcal{K}}$ and monotonically converges in $L_1$-norm to $\mathcal{I}_{\mathcal{K}}$ as its degree $d$ increases \cite{Ind2}. Hence, using Eq \eqref{eq: expect}, $E[\mathcal{P}_{\mathcal{K}}(x,q)]$ is an upper bound of $\mathbf{P_{risk}^*}$ and converges monotonically \cite{Ind1}. Similarly, $1-\mathcal{P}_{\bar{\mathcal{K}}}(x,q) \in \mathbb{R}_d[x,q] $ is a lower bound of $\mathcal{I}_{\mathcal{K}}$ and as its degree $d$ increases, $1-E[\mathcal{P}_{\bar{\mathcal{K}}}(x,q)]$ converges monotonically to $\mathbf{P_{risk}^*}$. $\blacksquare$\\

Let $\textbf{c}$ and $\bar{\textbf{c}}$ be the coefficient vectors of polynomials $\mathcal{P}_{\mathcal{K}}(x,q)$ and $\mathcal{P}_{\bar{\mathcal{K}}}(x,q)$, respectively, and $m^{xq}_i$ be the $i-$th moment of probability distribution $\mu_x\mu_q$, then lower and upper bounds of the risk in \eqref{risk_bound} can be written in terms of the weighted sum of the moments as {\small$E[\mathcal{P}_{\mathcal{K}}(x,q)]=\sum_{i} c_im^{xq}_i$} and {\small$1- E[\mathcal{P}_{\bar{\mathcal{K}}}(x,q)]=1-\sum_{i} \bar{c}_im^{xq}_i$}.
Hence, one can obtain the coefficients $\textbf{c}$ and $\bar{\textbf{c}}$ by solving SDP \eqref{Ind} in the offline step and then calculate the probability bounds for given probability distributions of uncertainties in the online step.\\

Note that the problem in \eqref{Ind} is a multivariate SOS optimization that looks for a polynomial of order $d$ in $(n+m)$-variate polynomial space (e.g., $S_{{n+m},d}$ unknown coefficients). Therefore, as the dimension of the original problem increases, the optimization problem \eqref{Ind} becomes computationally intractable. To avoid this, we present a procedure that requires solving a univariate SOS optimization (e.g., $S_{1,d}$ unknown coefficients). 

\subsection{Modified Moment-SOS Based Risk Bounds}
\label{sec:z}
In this section, to solve the risk estimation problem in \eqref{intro_P1}, we provide a procedure that requires a univariate approximation of the indicator function. For this purpose, we first consider the unsafe set $\chi(q)$ involving one polynomial (e.g., $\ell=1$) and then extend the obtained results to the set $\chi(q)$ involving multiple polynomials.

\subsubsection{Unsafe Set Involving One Polynomial}

Consider the given set $\chi$ in $\eqref{intro_set1}$ where $\ell=1$, i.e., $\chi(q):= \left\lbrace  x\in \reals^n: l_1 \leq \mathcal{P}(x,q)\leq u_1 \ \right\rbrace$. We define random variable $z \in \reals$ in terms of the polynomial of the set $\chi$ as 
\begin{equation}\label{Rand1}
z = \mathcal{P}(x,q)
\end{equation}
Random variable $z$ is a continuous function of the random variables $x$ and $q$; Therefore, its moments can be obtained in terms of the moments of probability distributions $\mu_x$ and $\mu_q$ as follows:
\begin{equation}\label{mom1}
m^z_{\alpha}=E[z^{\alpha}]=E[\mathcal{P}^{\alpha}(x,q)]=\sum_{i,j} a_{ij}m^x_im^q_{j}
\end{equation}where $m^z_{i},m^x_{i}$, and $m^q_{i}$ are the $i$-th moments of random variables $z$, $x$, and $q$, respectively, and $a_{ij}$ are the coefficients of polynomial $\mathcal{P}^{\alpha}(x,q)$.
Defining random variable $z$, the risk in \eqref{intro_P1} can be stated as:
\begin{align}
\mathbf{P_{risk}^*} := \mbox{Probability}_{\mu_z} \{ l_1 \leq z \leq u_1 \}
\end{align}where $\mu_z$ is the probability distribution of $z$. Note that, based on assumptions 2 and 3, random variable $z$ is supported on {$\mathcal{B}_z=[-1,1]$}. \\

According to \textit{Lemma} \ref{Lemma2}, the following results hold:
	\begin{equation} \label{risk_bound2}
	1- E[\mathcal{P}_{\bar{\mathcal{K}}}(z)]=1-\sum_{i=0}^{d} \bar{c}_{z_i}m^{z}_{i} \leq \mathbf{P_{risk}^*} \leq E[\mathcal{P}_{\mathcal{K}}(z)]=\sum_{i=0}^{d} c_{z_i}m^{z}_{i}
	\end{equation}
where {\small$\mathcal{K}=[l_1, u_1]$} and {\small $\bar{\mathcal{K}}=\mathcal{B}_z \setminus [l_1, u_1]$}. 
Univariate polynomials {\small$\mathcal{P}_{\mathcal{K}}(z) \in \reals_d[z]$} with coefficients {\small $c_{z_i}, \ i=0,...,d$} and {\small$\mathcal{P}_{\bar{\mathcal{K}}}(z)\in \reals_d[z] $} with coefficients {\small $\bar{c}_{z_i}, \ i=0,...,d$} are polynomial approximations of the indicator functions of the sets {\small$\mathcal{K}$} and {\small$\bar{\mathcal{K}}$}, respectively, that are obtained by solving a convex optimization problem similar to \eqref{Ind}, i.e., 
\begin{align}\label{Ind2}
& \mbox{min}_{\mathcal{P}_{\mathcal{K}}(z)\in \mathbb{R}_d[z]}  \int_{\mathcal{B}_z} \mathcal{P}_{\mathcal{K}}(z)dz \\
\mbox{s.t} &  \ \ \mathcal{P}_{\mathcal{K}}(z) \geq 1 \ \ \mbox{on} \ \ \mathcal{K} \subeqn  \\
& \ \ \mathcal{P}_{\mathcal{K}}(z) \geq 0 \ \ \mbox{on} \ \ \mathcal{B}_z \subeqn 
\end{align}
Note that the optimization problem in \eqref{Ind2} is a univariate SOS optimization. 
To improve the obtained risk bounds in \eqref{risk_bound2}, we use the Chebyshev polynomial basis instead of the standard basis to solve the optimization problem in \eqref{Ind2} and to represent the moments (\cite{Ind1,Ashk}). Hence, the new risk bounds read as:

	\begin{equation} \label{risk_bound3}
	1- E[\mathcal{P}_{T_{\bar{\mathcal{K}}}}]=1-\sum_{i=0}^{d} \bar{c}_{T_i}m^{z}_{T_i} \leq \mathbf{P_{risk}^*} \leq E[\mathcal{P}_{T_{\mathcal{K}}}]=\sum_{i=0}^{d} c_{T_i}m^{z}_{T_i}
	\end{equation}
where $\mathcal{P}_{T_{\mathcal{K}}}$ with coefficients $c_{T_i}, \ i=0,...,d$ and $\mathcal{P}_{T_{\bar{\mathcal{K}}}}$ with coefficients $\bar{c}_{T_i}, \ i=0,...,d$ are Chebyshev based polynomial approximations of the indicator functions of the sets $\mathcal{K}$ and $\mathcal{\bar{K}}$, respectively, and $m^{z}_{T_i}$ is the $i$-th moment of $z$ written in the Chebyshev basis.\\

\textbf{Illustrative Example:} Let $x$ be a single random variable with uniform probability distribution $\mu_x = U[-0.5, 0.5]$ representing the location of a ball. There is a moving hole $h=[q-0.8, \ q]$ where $q$ is a random variable with $\mu_q=Beta(3-\sqrt{2},3+\sqrt{2})$ probability distribution. 
We are interested in finding the probability that the ball lands in the hole, i.e., $\mbox{Probability}_{U,Beta} \{ x \in h\}$. Hence, the risk is defined as {\small $ \mathbf{P_{risk}^*} = \int_{-0.4\leq 0.5(x-q) \leq 0 } \mu_x\mu_qdxdq$.} We define a random variable as $z=0.5(x-q)$ and the sets $\mathcal{K}=[-0.4, 0]$ and $\mathcal{B}_z=[-1, 1]$. Then, {\footnotesize $\mathbf{P_{risk}^*} = \int_{\mathcal{K}} \mu_zdz= \int \mathcal{I}_{\mathcal{K}} \mu_zdz$}. The following result holds: {\footnotesize $ 1-\int_{-1}^{1} \mathcal{P}_{\bar{\mathcal{K}}}(z)\mu_z dz  \leq \mathbf{P_{risk}^*} \leq \int_{-1}^{1} \mathcal{P}_{\mathcal{K}}(z)\mu_z dz $} where {\small$\mathcal{P}_{\mathcal{K}}(z)$} and {\small$\mathcal{P}_{\bar{\mathcal{K}}}(z)$} are polynomial approximations of the indicator functions of the sets {\small$\mathcal{K}$} and {\small$\bar{\mathcal{K}}$}.\\

The $\alpha$-th moment of $z$ can be written in terms of the known moments of $x$ and $q$ as follows:
{\footnotesize $m^z_{\alpha}=E[z^{\alpha}]  =  E[(\frac{1}{2}(x-q))^{\alpha}]= \sum_{i=0}^{{\alpha}} \binom{\alpha}{i} (-1)^{\alpha-i}(\frac{1}{2})^{\alpha} m^x_i m^q  _{\alpha-i}$}, where the $i$-th moments of $x$ and $q$ are {\small$m^x_i=\frac{0.5^{i+1}-0.5^{i+1}}{i+1}$} , {\small $m^q_i=\frac{3-\sqrt{2}+i-1}{6+i-1}m^q_{i-1}$}, respectively. For example, the first 3 moments of $z$ described in the standard basis read as: {\small $m^z_0=1$, $m^z_1=0.5m^x_1-0.5m^q_1$, $m^z_2=	0.25m^x_2-0.5m^x_1m^q_1+0.25m^q_2$}. Also, the moments in the Chebyshev basis are
{\small $m^z_{T_0}=m^z_{0}$, $m^z_{T_1}=m^z_{1}$, $m^z_{T_2}=-m^z_{0}+2m^z_{2}$}.\\
	
\begin{figure}[t!]
	\centering
	\includegraphics[width=10cm, height=4cm]{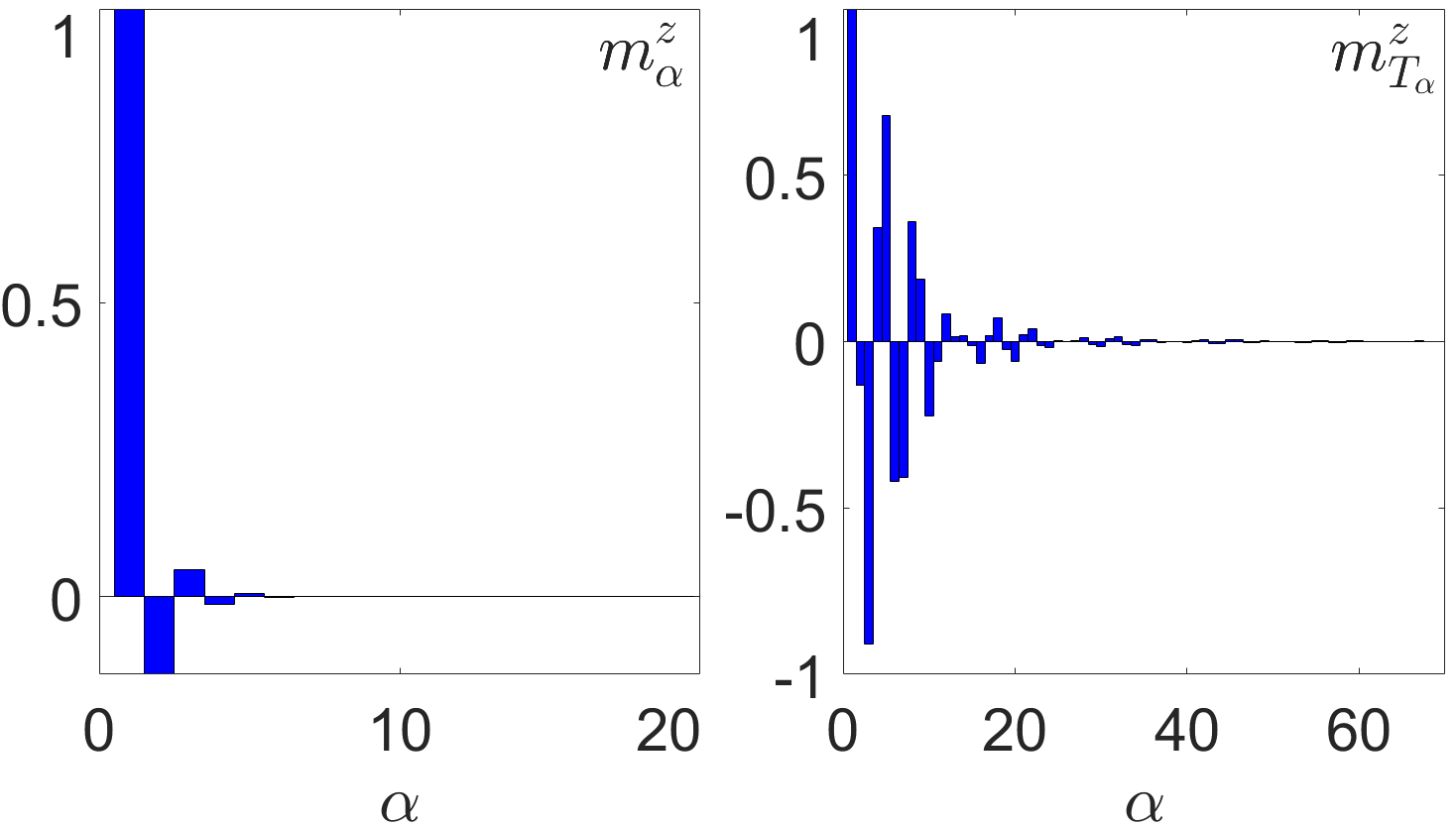}
	\caption{\small{Moments of random variable $z$ in the standard ($m^z_{\alpha}$ ) and Chebyshev basis ($m^z_{T_{\alpha}}$)}}
	\label{fig:E0m}
\end{figure}
\begin{figure}[t!]
	\centering
	\includegraphics[width=10cm, height=4cm]{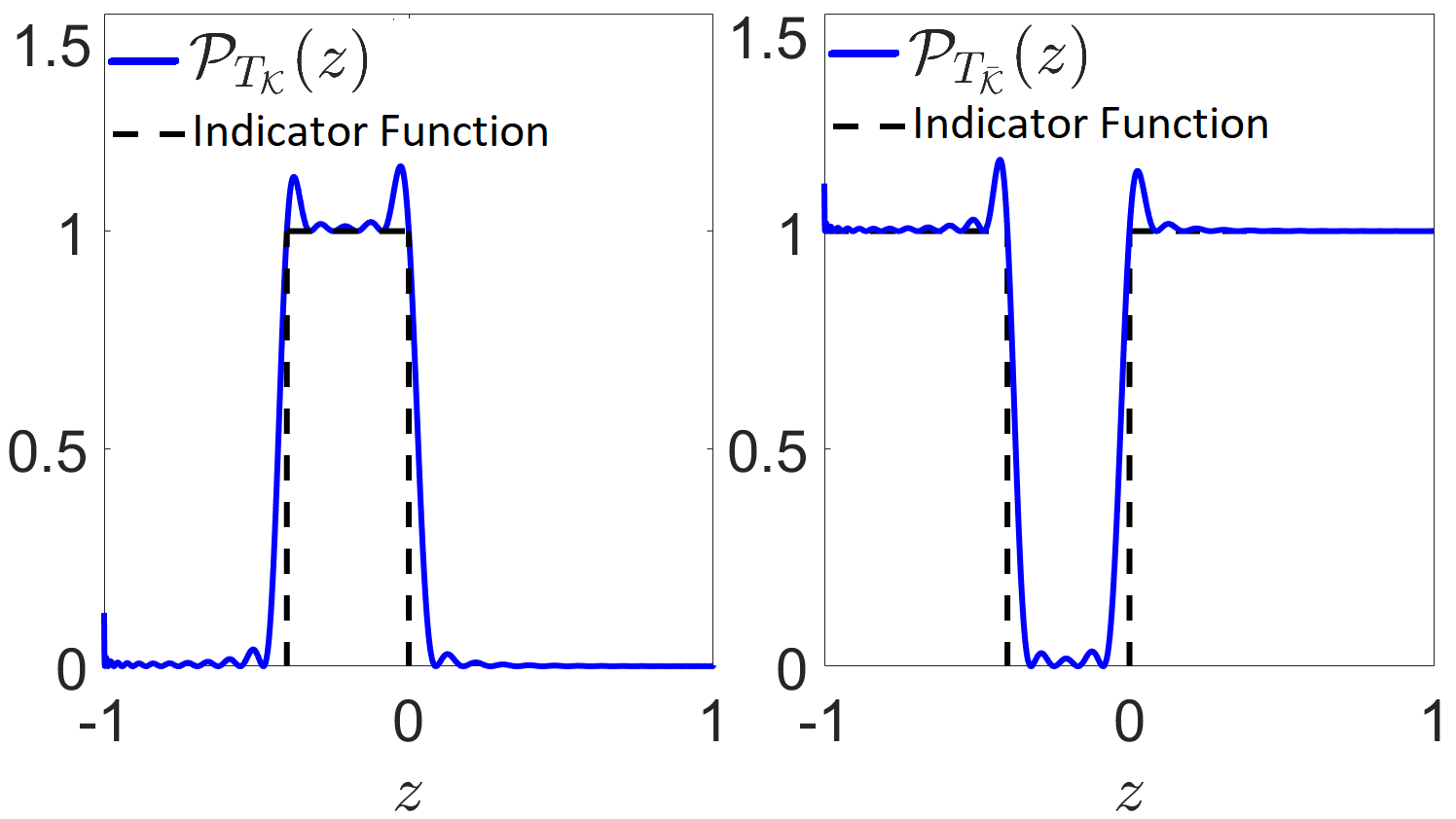}
	\caption{ {\small $\mathcal{P}_{T_{\mathcal{K}}}(z)$ and $\mathcal{P}_{T_{\bar{\mathcal{K}}}}(z)$, Chebyshev based polynomial approximations of the indicator functions of the sets $\mathcal{K}$, $\bar{\mathcal{K}}$, respectively.}} 
	\label{fig:E0}
\end{figure}
\begin{figure}[t!]
	\centering
	\includegraphics[width=10cm, height=4cm]{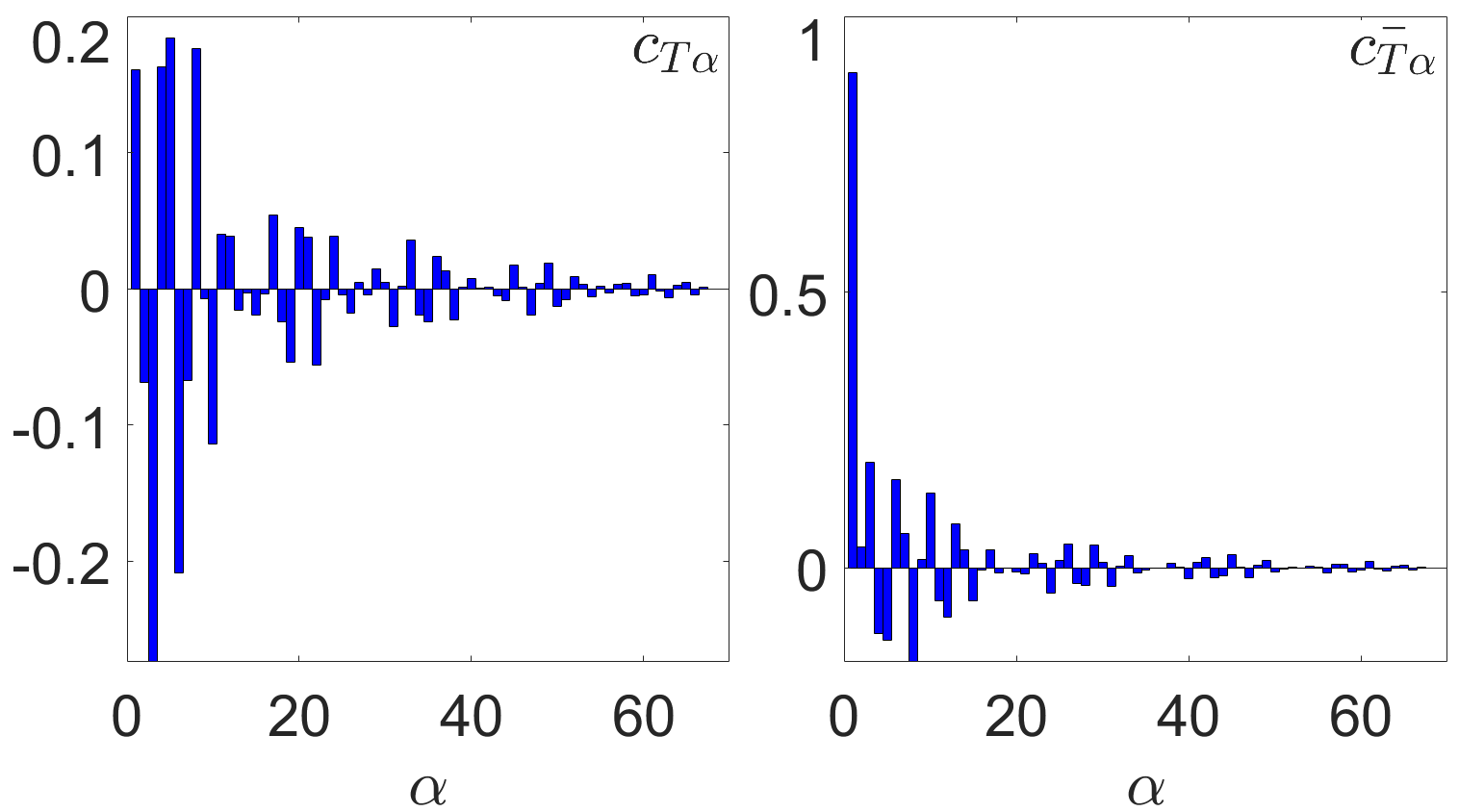}
	\caption{{\small $c_{T\alpha}$: coefficients of $\mathcal{P}_{T_{\mathcal{K}}}(z)$, $\bar{c_{T\alpha}}$: coefficients of  $\mathcal{P}_{T_{\bar{\mathcal{K}}}}(z)$}}
	\label{fig:E03}
\end{figure}

Figure \ref{fig:E0m} shows the moments of $z$ in the standard and Chebyshev basis up to the order $\alpha=66$. We solve optimization problem \eqref{Ind2} for the sets $\mathcal{K}$ and $\bar{\mathcal{K}}$ with $d=66$. The obtained polynomial approximations of the indicator functions and their coefficients are shown in Figures \ref{fig:E0} and \ref{fig:E03}, respectively. 
According to Eq \eqref{risk_bound3} the risk bounds are $[0.591, 0.798]$ while the true risk, approximated by the Monte-Carlo sampling method, is $0.7$. 
Table \ref{table:1} shows the obtained lower and upper bounds on the risk denoted by $p_l$ and $p_u$, respectively, for different polynomial degree $d$. 
According to \textit{Lemma} \ref{Lemma2} as $d$ increases, the obtained bounds converge to the true risk.  \\



It is shown that the Chebyshev basis could improve the resluts of SDPs (\cite{Ind1,Ashk}). In the risk estimation problem, the Chebyshev representation improves the results because i) Chebyshev based polynomial approximation of the indicator function reduces the oscillations on the boundary of the given set \cite{Ind1} e.g., points $x=-0.4$ and $x=0$ (Figure \ref{fig:E0}), ii) Chebyshev based representation of the moments affect the risk bounds more efficiently, e.g., moments in the standard basis $m^z_{\alpha}$ vanishes rapidly (Figure \ref{fig:E0m}).







%

%
\begin{table}[h!]
		\begin{center}
			\begin{tabular}{ |c|c|c|c|c|c|c| } 
				\hline
				d & 20 & 30 & 40 & 50 & 60 & 66   \\
				\hline
				${p_u}$  & 0.92 & 0.879 & 0.859 & 0.822 & 0.804 & 0.798 \\ \hline
				${p_l}$ & 0.401 & 0.485 & 0.511 & 0.562 & 0.586 & 0.591 \\ 
				\hline
			\end{tabular}
		\end{center}
		\caption{\small{Upper and lower bounds of the risk}}
		\label{table:1}
\end{table}

\subsubsection{Unsafe Set Involving Multiple Polynomials}
Consider the given set $\chi$ in $\eqref{intro_set1}$. We define random vector $Z \in \reals^{\ell}$ in terms of the polynomials of the set $\chi$ as: 
\begin{equation}\label{rand2}
Z=[z_1,...,z_\ell], \ \ z_j = \mathcal{P}_j(x,q), \ j=1,...,\ell
\end{equation}
The moment of order {\small$\overrightarrow{\alpha}=(\alpha_1,...,\alpha_{\ell})$} of $Z$ can be obtained in terms of the moments of probability distributions $\mu_x$ and $\mu_q$ as follows:
	\begin{equation}\label{mom22}
	m^{Z}_{\alpha_1,...,\alpha_{\ell}} =E[\prod_{j=1}^{\ell}z_j^{\alpha_j}]=E[\prod_{j=1}^{\ell}\mathcal{P}_j^{\alpha_j}(x,q)]
	=\sum_{i,j} a_{ij}m^x_im^q_{j}
	\end{equation}
where $m^x_{i}$ and $m^q_{i}$ are the $i$-th order moments of random variables $x$ and $q$, respectively, and $a_{ij}$ are the coefficients of polynomial {\small $\prod_{j=1}^{\ell}\mathcal{P}_j^{\alpha_j}(x,q)$}. Defining random vector $Z$, risk in \eqref{intro_P1} reads as
\begin{equation}
\mathbf{P_{risk}^*} := \mbox{Probability}_{\mu_Z} \{Z \in [l_1,u_1]\times[l_2,u_2]...\times[l_{{\ell}},u_{\ell}]\}
\end{equation}
where $\mu_Z$ is the probability distribution of $Z$. Hence, according to \textit{Lemma} \ref{Lemma2} the following result holds:
	\begin{equation} \label{risk_bound4}
	\mathbf{P_{risk}^*} \leq E[\prod_{j=1}^{\ell}\mathcal{P}_{\mathcal{K}_j}(z_j)]
	\end{equation}
where {\small$\mathcal{K}_j=[l_j, u_j]$, $j=1,...,\ell$} and univariate polynomials {\small$\mathcal{P}_{\mathcal{K}_j}(z_j)$}, {\small$j=1,...,\ell$} are the polynomial approximations of the indicator functions of the sets {\small$\mathcal{K}_j$, $j=1,...,\ell$} that are obtained by solving convex optimization problem \eqref{Ind2}. 
Note that {\small $\prod_j^{\ell}\mathcal{P}_{\mathcal{K}_j}(z_j)$} represents the polynomial approximation of the indicator function of the set {\small $[l_1,u_1]\times[l_2,u_2]...\times[l_{{\ell}},u_{\ell}]$}.\\

Let $ \mathcal{P}_{T_{{\mathcal{K}}_j}}$, $j=1,...,\ell$ denote the Chebyshev based polynomial approximations of the indicator functions of the sets {\small$\mathcal{K}_j$, $j=1,...,\ell$} and $m^{z}_{T_{i_1,...,i_{\ell}}}$ be the moment of order {\small$\overrightarrow{i}=(i_1,...,i_{\ell})$} of $Z$ written in the Chebyshev basis. Then, risk bound in \eqref{risk_bound4} reads as {\footnotesize  $\mathbf{P_{risk}^*} \leq \sum_{{i_1,...,i_{\ell}}} c_{T_{i_1,...,i_{\ell}}}m^{Z}_{T_{i_1,...,i_{\ell}}}$} where $c_{T_{i_1,...,i_{\ell}}}$ are the coefficients of polynomial {\small $\prod_{j=1}^{\ell} \mathcal{P}_{T_{{\mathcal{K}}_j}}$}.



\section{Implementation and Numerical Results}\label{sec:Res}
In this section, two numerical examples are presented that illustrate the performance of the proposed approach. We solve the SDP in \eqref{Ind2} to find the coefficients of Chebyshev based polynomial approximations of the indicator functions of the sets $\mathcal{K}$ and $\mathcal{\bar{K}}$. Also, we obtain the coefficients vector in \eqref{mom1} that maps the moments of uncertainties to the moments of random variable $z$. Obtaining these coefficients in the offline step, we calculate the risk bounds for any given uncertainties in real-time. Note that calculation of the risk bounds only requires multiplying the moment vector of uncertainties by the coefficient vectors calculated in the offline step. Hence, in the presence of time varying or state-dependent uncertainties, risk bounds can be updated by updating the moment information of uncertainties. For example, locations of dynamic obstacles can be modeled as time varying probabilistic uncertainties and the probability of collision at each time can be calculated by updating the moment information. \\

The computations in this section were performed on a computer with Intel i7 2.9GHz processors and 8 GB RAM. We use the Chebfun package \cite{Chebf} to work with univariate Chebyshev polynomials and also SeDuMi to solve the SDP in \eqref{Ind2}. We compare the proposed method with the moment-SDP based approach in \cite{Ind1} where, in the dual space, one needs to solve SOS optimization \eqref{Ind}. For this, we use GloptiPoly \cite{Glopti}, which is a MATLAB-based toolbox for moment-based SDP, and Mosek SDP solver. In all the tables, $d$ denotes the degree of the polynomial approximation of the indicator function, $p_u$ and $p_l$ denote upper and lower bounds on the risk, respectively, $t_u$ and $t_l$ denote computation time in seconds required for computing $p_u$ and $p_l$, respectively.\\

In this paper, we assume that semialgebraic representations of the obstacles are given. One can use the SOS based approaches in (\cite{obst1,Ind2}), to construct semialgebraic representations of obstacles from point cloud data obtained by sensors. In this case, additional constraints on the polynomials should be added to satisfy \textit{Assumption} 3. Note that defined random variable $z$ in the Section \ref{sec:z} is supported in $[-1,1]$; Hence, it's moment sequence is bounded in $[-1,1]$. We note, however, that describing the high order moments of $z$ in the Chebyshev basis could become numerically unstable. 
This is due to the large coefficients of the linear map between the Chebyshev and the standard basis that results in finite-precision floating-point error. 
Fixing this issue requires an appropriate rescaling of the Chebyshev basis \cite{Res}. 
In the provided numerical examples, we use the polynomial degree $d$ that results in bounded moments in the Chebyshev basis i.e., $-1\leq m^z_{T_{\alpha}}\leq 1, \alpha=0,...,d$.\\
\\

{\textbf{A. Example 1:}} Consider the motivating example in section \ref{Formu}. Based on the proposed approach, we need to find the following probability in terms of random variable $z$: $\mathbf{P_{risk}^*} = 
\mbox{Probability}_{\mu_z} \{-0.1\leq z \leq 0.2\}$. The polynomials $\mathcal{P}_{T_{\mathcal{K}}}$ and $\mathcal{P}_{T_{\bar{\mathcal{K}}}}$ for the set $\mathcal{K}=[-0.1, 0.2]$ are obtained by solving the optimization problem in \eqref{Ind2} with $d=88$. The moments of $z$ are obtained in terms of known moments of $x_{r1},x_{r2}$, and $q$ in the Chebyshev basis. Obtaining $m^z_{T_i}$ and coefficients $c_{T_i}$ and $\bar{c}_{T_i}$, the risk bounds are computed using \eqref{risk_bound3} as $[0.169, 0.48]$ while the true risk, approximated by the Monte-Carlo sampling method, is $0.32$. The computation time required to solve the univariate SOS optimization and calculate the risk in the online step are less than $\approx 17 (s)$ and $0.1(s)$, respectively.
Also, implementing the approach in \cite{Ind1}, the obtained upper and lower bounds of the risk and the computation time for different polynomial orders $d$ are reported in Table \ref{tab:ex1}. Note that one needs to repeat the heavy computations of the multivariate SOS in \cite{Ind1}, each time that moments of uncertainties change. However, in our approach one needs to only update the moment information and use the previously calculated coefficients $c_{T_i}$ and $\bar{c}_{T_i}$ to update the risk bounds. Proposed univariate SOS achieves better risk bounds in much less computation time.\\

\begin{table}
	\renewcommand{\arraystretch}{1.2}
		\begin{center}
			\begin{tabular}{| l | l | l | l |}
				\hline
					\multicolumn{4}{|c|}{{multivariate SOS }} \\ \thickhline
				$d$ &  10 & 20 & 30 \\ \hline	
				$p_u$& 0.54  & 0.50 & 0.495\\ \hline
				$t_u(s)$ & $\approx$2.6  & $\approx$76  & $\approx$3689 \\ \hline
				$p_l$ & 0.13 & 0.15 &0.161 \\ \hline
				$t_l(s)$ & $\approx$4.5  & $\approx$70  & $\approx$3156   \\ \thickhline
			\end{tabular}
					\begin{tabular}{| l | l |}
				\hline
				\multicolumn{2}{|c|}{{proposed univariate SOS}} \\ \thickhline
				$d$ &  88 \\ \hline	
				$p_u$& 0.48 \\ \hline
				$t_u(s)$ & $\approx$17 \\ \hline
				$p_l$ & 0.169 \\ \hline
				$t_l(s)$ & $\approx$15  \\ \thickhline	
			\end{tabular}
		\end{center}
	\vspace{3 mm}
	\caption{{\small Results of proposed univariate SOS and multivariate SOS in \cite{Ind1} for Example 1}}
	\label{tab:ex1}
\end{table}


\begin{figure}[t!]
	\centering
	\includegraphics[width=10cm, height=4cm]{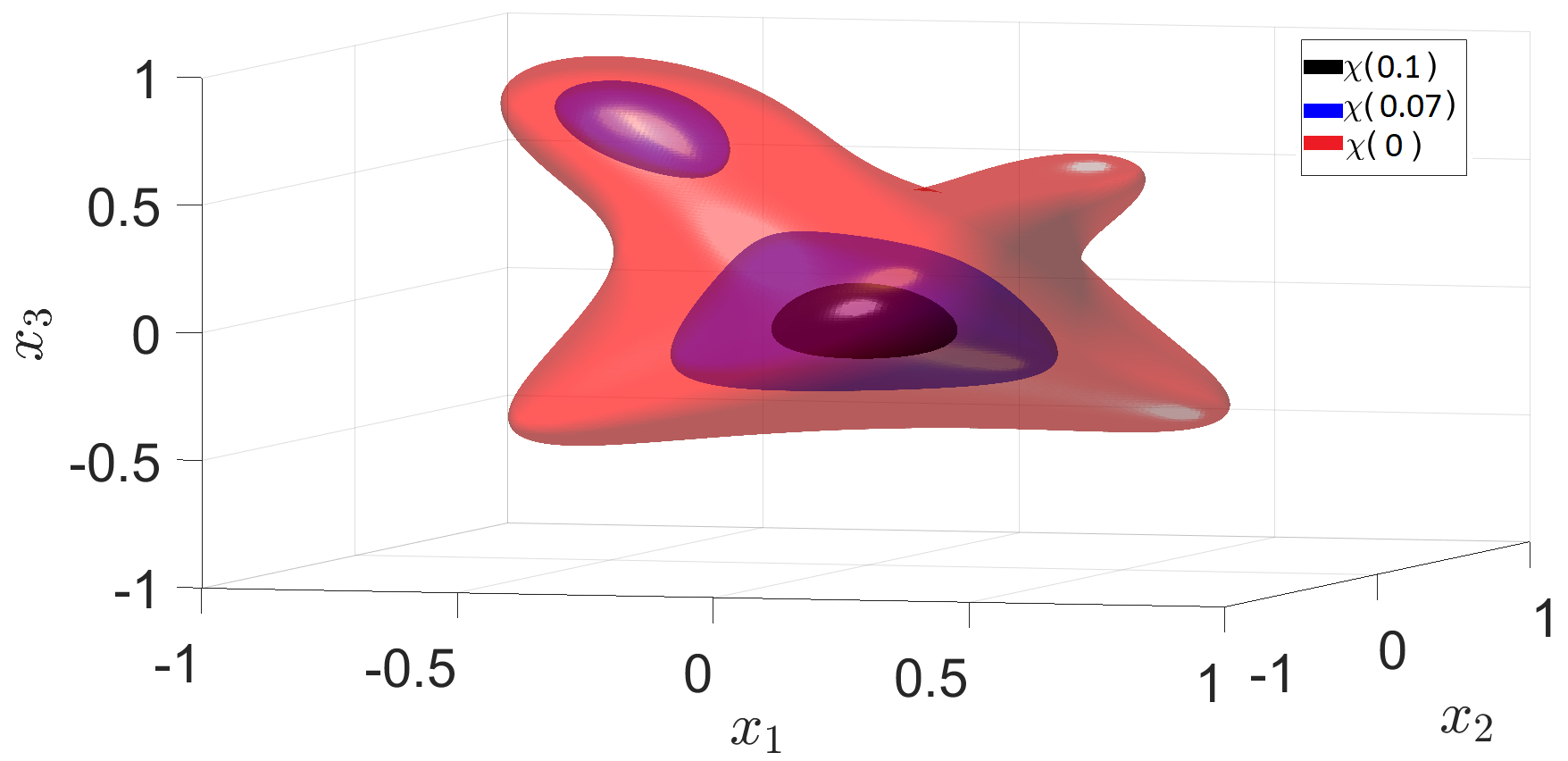}
	\caption{{\small Non-convex uncertain obstacle for $q=0$ (red), $q=0.07$ (blue), and $q=0.1$ (black)}}
	\label{fig:E3-1}
\end{figure}

{\textbf{B. Example 2:}} a non-convex uncertain set, shown in Figure \ref{fig:E3-1}, is described as $\chi(q)=\{x\in\reals^3:  0.84 \leq              
\mathcal{P}(x_1,x_2,x_3,q)
\leq 1
\}$, where 
{\scriptsize $\mathcal{P}=0.9487722614-0.0022x_1-0.0042x_2-0.0457x_3-0.3877x_1^2+0.0405x_1x_2-0.3105x_2^2-0.0537x_1x_3-0.0179x_2x_3-0.4094x_3^2-0.1059x_1^3-0.0212x_1^2x_2+0.0906x_1x_2^2-0.0543x_2^3+0.1451x_1^2x_3-1.8302x_1x_2x_3+0.1135x_2^2x_3-0.1096x_1x_3^2+0.1205x_2x_3^2+0.3407x_3^3-0.3285x_1^4-0.1338x_1^3x_2+0.4847x_1^2x_2^2+0.1127x_1x_2^3-0.3495x_2^4+0.0394x_1^3x_3+0.0149x_1^2x_2x_3-0.0051x_1x_2^2x_3-0.0594x_2^3x_3+0.5418x_1^2x_3^2-0.0659x_1x_2x_3^2+0.4840x_2^2x_3^2+0.0085x_1x_3^3+0.0657x_2x_3^3-0.3076x_3^4+0.1268x_1^5+0.0058x_1^4x_2-0.1012x_1^3x_2^2+0.0070x_1^2x_2^3+0.0053x_1x_2^4+0.0718x_2^5-0.0226x_1^4x_3+0.7338x_1^3x_2x_3-0.0716x_1^2x_2^2x_3+0.7226x_1x_2^3x_3-0.2075x_2^4x_3+0.0378x_1^3x_3^2-0.0139x_1^2x_2x_3^2+0.0224x_1x_2^2x_3^2-0.0566x_2^3x_3^2-0.0773x_1^2x_3^3+0.7345x_1x_2x_3^3+0.0955x_2^2x_3^3+0.0399x_1x_3^4-0.0653x_2x_3^4-0.3173x_3^5-q$} 
and $q$ is an uncertain parameter with uniform probability distribution on $[0,0.1]$. Also, 
$x_{1}$, $x_{2}$, and $x_{3}$ have uniform probability distributions on $[-0.4,0.4]$. We want to find the risk defined as $\mathbf{P_{risk}^*} = \mbox{Probability}_{\mu_{x_1},\mu_{x_2},\mu_{x_3},\mu_{q}} \{ x \in \chi(q) \}$. Based on the provided approach, we need to find the following probability in terms of new random variable $z$: $\mathbf{P_{risk}^*} = 
\mbox{Probability}_{\mu_z} \{0.84\leq z \leq 1\}$. The moments of $z$ are obtained in terms of known moments of $x_{1},x_{2},x_{3}$, and $q$ in the Chebyshev basis. The polynomials $\mathcal{P}_{T_{\mathcal{K}}}$ and $\mathcal{P}_{T_{\bar{\mathcal{K}}}}$ for the sets $\mathcal{K}=[0.84, 1]$ and $\bar{\mathcal{K}}=[-1,1]\setminus [0.84, 1]$ are obtained by solving optimization problem \eqref{Ind2} with $d=48$. Obtaining the moments $m^z_{T_i}$ and coefficients $c_{T_i}$ and $\bar{c}_{T_i}$, the risk bounds are obtained using \eqref{risk_bound3} as $[ 0.25, 0.77]$ while the true risk, approximated by the Monte-Carlo sampling method, is $0.519$. The computation time required to solve the univariate SOS optimization and calculate the risk in the online step are less than $\approx 5 (s)$ and $0.1(s)$, respectively. Also, implementing the approach in \cite{Ind1}, the obtained upper and lower bounds of the risk and the computation time for different polynomial orders $d$ are reported in Table \ref{tab:ex2}. For $d=30$, we receive an "out of memory" error due to the large size of the SDP.

\begin{table}
	\renewcommand{\arraystretch}{1.2}
		\begin{center}
			\begin{tabular}{| l | l | l | l |}
				\hline
					\multicolumn{4}{|c|}{{multivariate SOS }} \\ \thickhline
				$d$ &  10 & 20 & 30 \\ \hline	
				$p_u$& 0.81  & 0.78 & --\\ \hline
				$t_u(s)$ & $\approx$12  & $\approx$7459  & -- \\ \hline
				$p_l$ & 0.189 & 0.239 &-- \\ \hline
				$t_l(s)$ & $\approx$11  & $\approx$6657  & --   \\ \thickhline
			\end{tabular}
					\begin{tabular}{| l | l |}
				\hline
				\multicolumn{2}{|c|}{{proposed univariate SOS}} \\ \thickhline
				$d$ &  48 \\ \hline	
				$p_u$& 0.77 \\ \hline
				$t_u(s)$ & $\approx$5 \\ \hline
				$p_l$ & 0.25 \\ \hline
				$t_l(s)$ & $\approx$5  \\ \thickhline	
			\end{tabular}
		\end{center}
	\vspace{3 mm}
	\caption{{\small Results of proposed univariate SOS and multivariate SOS in \cite{Ind1} for Example 2}}
	\label{tab:ex2}
\end{table}

\section{Conclusion}\label{sec:Con}

In this paper, we consider the probability estimation of the safety constraints violation in the presence of bounded uncertainties with arbitrary probability distributions. Safety constraints are represented by a non-convex set defined by polynomial inequalities. To solve this problem, we use a moment-based representation of probability distributions. Upper and lower bounds of the risk are computed as a weighted sum of the moments of the probability distributions of uncertainties. The weights are obtained in the offline step by solving a univariate sum of squares optimization problem in the Chebyshev basis. Numerical examples on probabilistic collision checking problem in uncertain environments are provided that show the performance of the proposed method. For the future work, we will use the proposed method in probabilistic motion planning to evaluate the risk of the designed maneuvers for robots.


{}

\end{document}

%% file: defs.tex
\def\cK{\mathcal{K}}

\def\cP{\mathcal{P}}

\def\smskip{\smallskip}

\def\texitem#1{\par\smskip\noindent\hangindent 25pt
               \hbox to 25pt {\hss #1 ~}\ignorespaces}

\def\norm#1{\|#1\|}

\newcommand{\BEAS}{\begin{eqnarray*}}
\newcommand{\EEAS}{\end{eqnarray*}}
\newcommand{\BEA}{\begin{eqnarray}}
\newcommand{\EEA}{\end{eqnarray}}
\newcommand{\BEQ}{\begin{eqnarray}}
\newcommand{\EEQ}{\end{eqnarray}}
\newcommand{\BIT}{\begin{itemize}}
\newcommand{\EIT}{\end{itemize}}
\newcommand{\BNUM}{\begin{enumerate}}
\newcommand{\ENUM}{\end{enumerate}}

\newcommand{\BA}{\begin{array}}
\newcommand{\EA}{\end{array}}

\newcommand{\reals}{\mathbb{R}}

\newif\ifpagenumbering
\pagenumberingtrue

\pagenumberingfalse

\newsavebox{\theorembox}
\newsavebox{\lemmabox}
\newsavebox{\remarkbox}
\newsavebox{\assbox}
\savebox{\theorembox}{\noindent\bf Theorem}
\savebox{\lemmabox}{\noindent\bf Lemma}
\savebox{\remarkbox}{\noindent\bf Remark}
\savebox{\assbox}{\noindent\bf Assumption}

